\documentclass[11pt]{amsart}

\usepackage[usenames]{color}

\usepackage{amsmath,times,pslatex,enumitem}
\usepackage{amsthm}
\usepackage{amsfonts}
\usepackage{amssymb}
\usepackage{latexsym}

\usepackage{hyperref}

\newtheorem{theorem}{Theorem}
\newtheorem{corollary}{Corollary}

\newtheorem{problem}{Problem}

\newtheorem{proposition}{Proposition}

\newcommand{\Z}{\mathbb{Z}}

\newcommand{\Q}{\mathbb{Q}}

\newcommand{\F}{\mathbb{F}}

\newcommand{\R}{\mathbb{R}}

\begin{document}

\title{On two-generator subgroups in $SL_2(\Z)$, $SL_2(\Q)$, and  $SL_2(\R)$}

\author[]{Anastasiia Chorna}
\address{Department of Mathematics, The City  College  of New York, New York,
NY 10031} \email{anastasiyach.a@gmail.com}

\author[]{Katherine Geller}
\address{Department of Mathematics, The City  College  of New York, New York,
NY 10031} \email{geller.katherine@gmail.com}

\author[]{Vladimir Shpilrain}
\address{Department of Mathematics, The City  College  of New York, New York,
NY 10031} \email{shpil@groups.sci.ccny.cuny.edu}
\thanks{Research of Vladimir Shpilrain was partially supported by
the NSF grant CNS-1117675  and  by the ONR (Office of Naval
Research) grant N000141512164}

\begin{abstract}
We consider what some authors call ``parabolic M\"obius subgroups"
of matrices over $\Z$, $\Q$, and $\R$ and focus on the membership
problem in these subgroups and complexity of relevant algorithms.

\end{abstract}

\maketitle

\section{Introduction: two theorems of Sanov}

Denote  $A(k) = \left(
 \begin{array}{cc} 1 & k \\ 0 & 1 \end{array} \right) , \hskip .2cm B(k) = \left(
 \begin{array}{cc} 1 & 0 \\ k & 1 \end{array} \right).$ In an old
paper \cite{Sanov}, I. N. Sanov proved two simple yet remarkable
theorems:

\begin{theorem}\label{th1}
The subgroup of $SL_2(\Z)$ generated by $A(2)$ and $B(2)$ is free.
\end{theorem}

\begin{theorem}\label{th2}
The subgroup of $SL_2(\Z)$ generated by $A(2)$ and $B(2)$ consists
of {\it all} matrices of the form  $\left(
\begin{array}{cc} 1+4n_1 & 2n_2 \\ 2n_3 & 1+4n_4\end{array}
\right)$ with determinant 1, where all $n_i$ are arbitrary integers.
\end{theorem}

These two theorems together yield yet another proof of the fact that
the group $SL_2(\Z)$ is virtually free. This is because the group of
all invertible matrices of the form  $\left(
\begin{array}{cc} 1+4n_1 & 2n_2 \\ 2n_3 & 1+4n_4\end{array}
\right)$  obviously has finite index in $SL_2(\Z)$. Thus, we have:

\begin{corollary}\label{corollary1}
The group $SL_2(\Z)$ is virtually free.
\end{corollary}

There is another interesting corollary of Theorem \ref{th2}:

\begin{corollary}\label{corollary2}
The membership problem in the subgroup of $SL_2(\Z)$ generated by
$A(2)$ and $B(2)$ is solvable in constant time.
\end{corollary}

We note that this is, to the best of our knowledge, the only example
of a natural (and nontrivial) algorithmic problem in group theory
solvable in constant time. In fact, even problems solvable in
sublinear time are very rare, see \cite{Sublinear}, and in those
that are, one can typically get either ``yes" or ``no" answer in
sublinear time, but not both. Complexity of an input in our case is
the ``size" of a given matrix, i.e., the sum of the absolute values
of its entries. In light of Theorem \ref{th2}, deciding whether or
not a given matrix from $SL_2(\Z)$ belongs to the subgroup generated
by $A(2)$ and $B(2)$ boils down to looking at residues modulo 2 or 4
of the entries. The latter is decided by looking just at the last
one or two digits of each entry (assuming that the entries are given in the binary or, say, decimal form). We emphasize though that solving
this membership problem in constant time is only possible if an input
matrix is known to belong to $SL_2(\Z)$; otherwise one would have
to check that the determinant of a given matrix is equal to 1, which
cannot be  done in constant time, although can still be done in
sublinear time with respect to the complexity $|M|$ of an input
matrix $M$, as defined in the next Section \ref{results}, see
Corollary \ref{complexity}.

\section{Our results}
\label{results}

In this paper, we show that what would be a natural generalization
of Sanov's Theorem \ref{th2} to $A(k)$ and $B(k)$, $k \in \Z_+$, is
not valid for $k \ge 3$ and moreover, the subgroup generated by
$A(k)$ and $B(k)$ has infinite index in $SL_2(\Z)$ if $k \ge 3$.

\begin{theorem}\label{infinite}
The subgroup of $SL_2(\Z)$ generated by $A(k)$ and $B(k)$, $k \in
\Z, ~k \ge 3$, has infinite index in the group of all matrices of
the form $\left(
\begin{array}{cc} 1+k^2m_1 & km_2 \\ km_3 & 1+k^2m_4\end{array}
\right)$ with determinant 1.
\end{theorem}

The group of  all matrices of the above form, on the other hand,
obviously has finite index in $SL_2(\Z)$.

Our main technical result, proved in Section \ref{Peak}, is the
following

\begin{theorem} \label{peak}
Let $M = \left(
 \begin{array}{cc} m_{11} & m_{12} \\ m_{21} & m_{22} \end{array}
 \right)$ be a matrix from $SL_2(\R)$. Call ``elementary operations"
on $M$ the following 8 operations: multiplication of $M$ by either
$A(k)^{\pm 1}$ or by  $B(k)^{\pm 1}$, on the right or on the left.
\medskip

\noindent {\bf (a)} If $k \in \Z$ and $M$ belongs to the subgroup of
$SL_2(\Z)$ generated by $A(k)$ and $B(k)$, then it has the form
$\left(
\begin{array}{cc} 1+k^2n_1 & kn_2 \\ kn_3 & 1+k^2n_4\end{array}
\right)$ for some integers $n_i$.

If $k \in \R$ and  $M$ belongs to the subgroup of $SL_2(\R)$
generated by $A(k)$ and $B(k)$, then it has the form $\left(
\begin{array}{cc} 1+\sum_i k^in_i & \sum_j k^jn_j \\ \sum_r k^rn_r & 1+\sum_s k^sn_s\end{array}
\right)$ where all $n_i$ are integers and all exponents on $k$ are
positive integers.
\medskip

\noindent {\bf (b)} Let $k \in \R, ~k \ge 2$. If $M \in SL_2(\R)$
and there is a sequence of elementary operations that reduces
$\sum_{i,j}|m_{ij}|$, then there is a single elementary operation
that reduces $\sum_{i,j}|m_{ij}|$.

\noindent {\bf (c)} Let $k \in \Z, ~k \ge 2$. If $M \in SL_2(\Z)$
and no single elementary operation reduces $\sum_{i,j}|m_{ij}|$, then
either $M$ is the identity matrix or $M$ does not belong to the
subgroup generated by $A(k)$ and $B(k)$.

\end{theorem}

%
%

We also point out a result, similar to Theorem \ref{peak}, about the
{\it monoid} generated by $A(k)$ and $B(k)$ for $k>0$. Unlike
Theorem \ref{peak} itself, this result is trivial.

\begin{proposition}
Let $M = \left(
 \begin{array}{cc} m_{11} & m_{12} \\ m_{21} & m_{22} \end{array}
 \right)$ be a matrix from $SL_2(\R)$. Call ``elementary operations"
on $M$ the following 4 operations: multiplication of $M$ by either
$A(k)^{-1}$ or by  $B(k)^{-1}$, on the right or on the left.
\medskip

\noindent {\bf (a)} If $k \in \Z, ~k > 0,$ and $M$ belongs to the
{\it monoid} generated by $A(k)$ and $B(k)$, then it has the form
$\left(
\begin{array}{cc} 1+k^2n_1 & kn_2 \\ kn_3 & 1+k^2n_4\end{array}
\right)$ for some nonnegative integers $n_i$.

If $k \in \R, ~k > 0,$ and  $M$ belongs to the {\it monoid}
generated by $A(k)$ and $B(k)$, then it has the form $\left(
\begin{array}{cc} 1+\sum_i k^in_i & \sum_j k^jn_j \\ \sum_r k^rn_r & 1+\sum_s k^sn_s\end{array}
\right)$ where all $n_i$ are nonnegative  integers and all exponents
on $k$ are positive integers.
\medskip

\noindent {\bf (b)} Let $k \in \Z, ~k \ge 2$. If $M$ is a matrix
from $SL_2(\Z)$ with nonnegative entries and no elementary operation
reduces $\sum_{i,j} m_{ij}$, then either $M$ is the identity matrix
or  $M$ does not belong to the {\it monoid}  generated by $A(k)$ and
$B(k)$.

\end{proposition}

Thus, for example, the matrix $\left(
\begin{array}{cc} 5 & 4 \\ 6 & 5 \end{array}
\right)$ does not belong to the {\it monoid}  generated by $A(2)$
and $B(2)$, although it does belong to the {\it group} generated by
$A(2)$ and $B(2)$ by Sanov's Theorem \ref{th2}.


Theorem \ref{peak} yields a simple algorithm for the membership
problem in the subgroup generated by $A(k)$ and $B(k)$ in case $k
\in \Z, ~k \ge 2.$ We note in passing that in general, the subgroup
membership problem for $SL_2(\Q)$ is open, while in $SL_2(\Z)$ it is
solvable since $SL_2(\Z)$ is virtually free. The general solution,
based on the automatic structure of $SL_2(\Z)$ (see \cite{Epstein}),
is not so transparent and has quadratic time complexity (with respect to
the word length of an input). For our
special subgroups we have:

\begin{corollary}\label{complexity} Let $k \in \Z, ~k \ge 2$, and
let the complexity $|M|$ of a matrix $M = \left(
 \begin{array}{cc} m_{11} & m_{12} \\ m_{21} & m_{22} \end{array}
 \right)$ be the sum of all $|m_{ij}|$. There is an algorithm that
decides whether or not a given matrix $M \in SL_2(\Z)$ is in the
subgroup of $SL_2(\Z)$ generated by $A(k)$ and $B(k)$ (and if it
does, finds a presentation of $M$ as a group word in $A(k)$ and
$B(k)$) in time $O(n \cdot \log n)$, where  $n=|M|$.

\end{corollary}

\noindent {\bf Remark.} The relation between $|M|= \sum |m_{ij}|$
and the word length of $M$ (with respect to the standard generators
$A(1)$ and $B(1)$, say) is not at all obvious and is an  interesting
problem in its own right.

\medskip

Statement similar to Corollary \ref{complexity} holds also for the
{\it monoid} generated by $A(k)$ and $B(k)$, for any $k \in \Z, ~k
\ge 2$.

The $O(n \cdot \log n)$ is the worst-case complexity of the
algorithm referred to in Corollary \ref{complexity}. It would be
interesting to find out what the {\it generic-case complexity} (in
the sense of \cite{KMSS}) of this algorithm is. Proposition 1 in
\cite{BSV} tacitly suggests that this complexity might be, in fact,
sublinear in $n=|M|$, which would be a really interesting result, so
we ask:

\begin{problem}\label{generic} Is the generic-case complexity of the algorithm
claimed in Corollary \ref{complexity} sublinear in $|M|$?
\end{problem}

We note that, unlike the algorithms with low generic-case complexity
considered in \cite{KMSS}, this algorithm has a good chance to have
low generic-case complexity giving both ``yes" and ``no" answers,
see our Section \ref{Corollary} for more details.

Finally, we note that if $M$ is in the subgroup generated by $A(k)$
and $B(k)$, $k \ge 2$, then the presentation of $M$ as a group word
in $A(k)$ and $B(k)$ is unique since the group generated by $A(k)$
and $B(k)$ is known to be free for any $k \in \R, ~k \ge 2$, see
e.g. \cite{JP}. On the other hand, the group generated by $A(1)$ and
$B(1)$ (i.e., the whole group $SL_2(\Z)$) is not free. This implies,
in particular, that for any integer $n \ge 1$, the group generated
by $A(\frac{1}{n})$ and $B(\frac{1}{n})$ is not free because it
contains both matrices $A(1)$ and $B(1)$. Many examples of rational
$k, ~ 0 < k < 2$, for which the subgroup of $SL_2(\Q)$ generated by
$A(k)$ and $B(k)$ is not free were found over the years, starting
with \cite{L-U}; see a recent paper \cite{Gutan} for more
references. (We can single out the paper \cite{Beardon} where the
question of non-freeness for this subgroup was reduced to
solvability of particular Diophantine equations.) In particular, it
is known that for any $k, ~ 0 < k < 2$, of the form $\frac{m}{mn+1}$
or $\frac{m+n}{mn}, ~m, n \in Z_+$, the group generated by $A(k)$
and $B(k)$ is not free. This includes $k=\frac{2}{3}, \frac{3}{2},
\frac{3}{7}$, etc. Also, if the group is not free for some $k$, then
it is not free for any $\frac{k}{n}, ~n \in Z_+$.

The following problem, however, seems to be still open:

\begin{problem} (Yu. Merzlyakov \cite{Problems}, \cite{Kourovka})
For which rational $k, ~ 0 < k < 2$, is the group generated by
$A(k)$ and $B(k)$ free? More generally, for which algebraic $k, ~ 0
< k < 2$, is this group free?
\end{problem}

To the best of our knowledge, there are no known examples of a
rational $k, ~ 0 < k < 2$, such that the group generated by $A(k)$
and $B(k)$ is free. On the other hand, since any matrix from this
group has the form $ \left(
\begin{array}{cc} p_{11}(k) & p_{12}(k)\\
  p_{21}(k) & p_{22}(k) \end{array}
\right)$ for some polynomials $p_{ij}(k)$ with integer coefficients,
this group is obviously free if $k$ is transcendental. For the same
reason, if $r$ and $s$ are algebraic numbers that are Galois
conjugate over $\Q$, then the group generated by $A(r)$ and $B(r)$
is free if and only if the group generated by $A(s)$ and $B(s)$ is.
For example, if $r=2-\sqrt 2$, then $A(r)$ and $B(r)$ generate a
free group because this $r$ is Galois conjugate to $s=2+\sqrt 2 >
2$. More generally, $A(r)$ and $B(r)$ generate a free group for
$r=m- n\sqrt 2$, and therefore also for $r=k \cdot (m- n\sqrt 2)$,
with arbitrary positive $k, m, n \in \Z$. This implies, in
particular, that the set of algebraic $r$ for which the group is
free is dense in $\R$ because $(m- n\sqrt 2)$ can be arbitrarily
close to 0. All these $r$ are irrational though.

%
%

\section{Peak reduction}
\label{Peak}

Here we prove Theorem \ref{peak} from Section \ref{results}. The
method we use is called {\it peak reduction} and goes back to
Whitehead \cite{Wh}, see also \cite{L-S}. The idea is as follows.
Given an algorithmic problem that has to be solved, one first
somehow defines {\it complexity} of possible inputs. Another
ingredient is a collection of {\it elementary operations} that can
be applied to inputs. Thus, we now have an action of the semigroup
of elementary operations on the set of inputs. Usually, of
particular interest are elements of minimum complexity in any given
orbit under this action. The main problem typically is to find these
elements of minimum complexity. This is where the peak reduction
method can be helpful. A crucial observation is: if there is a
sequence of elementary operations (applied to a given input) such
that at some point in this sequence the complexity goes up (or
remains unchanged) before eventually going down, then there must be
a pair of {\it subsequent} elementary operations in this sequence (a
``peak") such that  one of them increases the complexity (or leaves
it unchanged), and then the other one decreases it. Then one tries
to prove that such a peak can always be reduced, i.e., if there is
such a pair, then there is also a single elementary operation that
reduces complexity. This will then imply that there is a ``greedy"
sequence of elementary operations, i.e., one that reduces complexity
{\it at every step}. This will yield an algorithm for finding an
element of minimum complexity in a given orbit.

In our situation, inputs are matrices from $SL_2(\R)$. For the
purposes of the proof of Theorem \ref{peak}, we define complexity of
a matrix $M = \left( \begin{array}{cc} m_{11} & m_{12} \\ m_{21} &
m_{22} \end{array} \right)$ to be the maximum of all $|m_{ij}|$.
Between two matrices with the same $\max |m_{ij}|$, the one with the
larger $\sum_{i,j}|m_{ij}|$ has higher complexity. We will see,
however, that in case of $2 \times 2$ matrices with determinant 1,
the ``greedy" sequence of elementary operations would be the same as
if we defined the complexity to be just $\sum_{i,j}|m_{ij}|$.

Elementary operations in our situation are multiplications of a
matrix by either $A(k)^{\pm 1}$ or by  $B(k)^{\pm 1}$, on the right
or on the left. They correspond to elementary row or column
operations; specifically, to operations of the form $(row_1 \pm  k
\cdot row_2)$, $(row_2 \pm  k \cdot row_1)$, $(column_1 \pm  k \cdot
column_2)$, and  $(column_2 \pm  k \cdot column_1)$.

%


We now get to

\noindent {\bf Proof of Theorem \ref{peak}.} Part (a) is established
by an obvious induction on the length of a group word representing a
given element of the subgroup generated by $A(k)$ and $B(k)$. Part
(c) follows from part (b). We omit the details and proceed to part
(b).

We are going to consider various pairs of subsequent elementary
operations of the following kind: the first operation increases the
maximum of $|m_{ij}|$ (or leaves it unchanged), and then the second
one reduces it. We are assuming that the second elementary operation
is not the inverse of the first one.

In each case like that, we show that either the maximum of
$|m_{ij}|$ in the given matrix could have been reduced by just a
single elementary operation (and then $\sum |m_{ij}|$ should be reduced, too, to keep the determinant unchangedk), or $\sum |m_{ij}|$ could have been
reduced by a single elementary operation leaving the maximum of
$|m_{ij}|$ unchanged. Because of a ``symmetry", it is sufficient to
consider the following cases.

First of all, we note that since the determinant of $M$ is equal to
1, there can be 0, 2, or 4  negative entries in $M$. If there are 2
negative entries, they can occur either in the same row, or in the
same column, or on the same diagonal. Because of the symmetry, we
only consider the case where two negative entries are in the first
column and the case where they are on the main diagonal. Also, cases
with 0 and 4 negative entries are symmetric, so we only consider the
case where there are no negative entries.

It is also convenient for us to single out the case where $M$ has
two zero entries, so we start with

\medskip

\noindent {\bf Case 0.} There are two zero entries in $M$. Since the
determinant of $M$ is 1, the two  nonzero entries should be on a
diagonal and their product should be $\pm 1$. If  they are not on
the main diagonal, then $M^4 = I$,  in which case $M$ cannot belong
to the subgroup generated by $A(k)$ and $B(k)$, $k \ge 2$, since
this subgroup is free.

Now suppose that the two  nonzero entries are on the main diagonal,
so $M = \left(
\begin{array}{cc} x & 0\\
 0 & \frac{1}{x} \end{array}
\right)$ for some $x \in \R, ~x \ne 1$. Without loss of generality,
assume $x > 0$. We are going to show, by way of contradiction, that
such a matrix is not in the subgroup generated by $A(k)$ and $B(k)$.
We have:  $MA(k)M^{-1} = \left(
\begin{array}{cc} 1 & x^2k\\
 0 & 1 \end{array}
\right),$ and for some $r \in \Z_+$ we have $A(k)^{-r}MA(k)M^{-1} =
\left(
\begin{array}{cc} 1 & yk\\
 0 & 1 \end{array}
\right),$ where $0 < yk \le k$. If  $yk = k$, then we have a
relation $A(k)^{-r}MA(k)M^{-1} = A(k)$, so again we have a
contradiction with the fact that the subgroup generated by $A(k)$
and $B(k)$ is free. Now let $0 < yk < k$, so $0 < y < 1$, and let $C
=  \left(\begin{array}{cc} 1 & yk\\
 0 & 1 \end{array}\right)$. We claim that $C$ does not belong to the
subgroup generated by $A(k)$ and $B(k)$. If it did, then so would
the matrix $T(m, n) = A(k)^{-m}C^n =  \left(\begin{array}{cc} 1 & (ny-m)k\\
 0 & 1 \end{array}\right)$ for any $m, n \in \Z_+$. Since $(ny-m)k$
 can be arbitrarily close to 0, the matrix $T(m, n)$ can be arbitrarily close to
the identity matrix, which contradicts the well-known fact (see e.g.
\cite{Jorgensen}) that the group generated by $A(k)$ and $B(k)$ is
{\it discrete} for any $k \ge 2$.



\medskip

In what follows, we assume that all matrices under consideration
have at most one zero entry. Even though we use strict inequalities
for all entries of a matrix, the reader should keep in mind that one
of the inequalities may  be not strict; this does not affect the
argument.

\smallskip

\noindent {\bf Case 1.} There are 2 negative entries, both in the
first column. Thus,  $m_{11}<0, m_{21}<0, m_{12}>0, m_{22}>0$.

\medskip

\noindent {\bf Case 1a.} Two subsequent elementary operations
reducing some entry after increasing it are both $(column_1 + k
\cdot column_2)$. If, after one operation $(column_1 + k \cdot
column_2)$, the element $m_{11}$, say, becomes positive, then
$m_{21}$ should become positive, too, for the determinant to remain
unchanged. Then, after applying $(column_1 + k \cdot column_2)$ one
more time, new $|m_{11}|$ will become greater than it was, contrary
to the assumption.

If, after one operation $(column_1 + k \cdot column_2)$, $m_{11}$
remains negative, then this operation reduces $|m_{11}|$, and this
same operation should also reduce  $|m_{21}|$  for the determinant
to remain unchanged. Indeed, the determinant is $m_{11}m_{22} -
m_{12}m_{21} = 1$. If $|m_{11}|$ decreases while $m_{11}$ remains
negative, then the value of $m_{11}m_{22}$ increases (but remains
negative). Therefore, the value of $m_{12}m_{21}$ should increase,
too, for the difference to remain unchanged. Since $m_{21}$ should
remain negative, this implies that $|m_{21}|$ should decrease, hence
$\sum |m_{ij}|$ decreases after one operation $(column_1 + k \cdot
column_2)$.

The same kind of argument works in the case where both operations
are  $(row_1 - k \cdot row_2)$,

If both operations are $(row_1 + k \cdot row_2)$, or $(column_1 - k
\cdot column_2)$, then no $|m_{ij}|$ can possibly decrease since
$m_{11}<0, m_{21}<0, m_{12}>0, m_{22}>0$.

\smallskip

\noindent {\bf Case 1b.} Two subsequent elementary operations are:
$(row_1 + k \cdot row_2)$, followed by $(column_1 + k \cdot
column_2)$. The result of applying these two operations to the
matrix $M$ is:  $ \left(
\begin{array}{cc} m_{11} + k \cdot m_{21} +k \cdot m_{12} + k^2 m_{22} & m_{12}+k \cdot m_{22}\\
  m_{21}+k \cdot m_{22} & m_{22} \end{array}
\right)$. This case is nontrivial only if the first operation
increases the absolute value of the  element in the top left corner
(or leaves it unchanged), and then the second operation reduces it.
But then the second operation should also reduce the absolute value
of $m_{21}$ for the determinant to remain unchanged. This means a
single operation  $(column_1 + k \cdot column_2)$ would reduce
$|m_{21}|$ to begin with,  and this same operation should also
reduce $|m_{11}|$ for the determinant to remain unchanged. Thus, a
single elementary operation would reduce the complexity of $M$.

The same argument takes care of any of the following pairs of
subsequent elementary operations: $(row_1 \pm k \cdot row_2)$,
followed by $(column_1 \pm k \cdot column_2)$, as well as  $(row_1
\pm k \cdot row_2)$, followed by $(column_2 \pm k \cdot column_1)$.

\smallskip

\noindent {\bf Case 1c.} Two subsequent elementary operations are:
 $(row_1 + k \cdot row_2)$, followed by $(row_2 + k \cdot row_1)$. The result of
applying these two operations to the matrix $M$ is:  $ \left(
\begin{array}{cc} m_{11}+k \cdot m_{21} & m_{12}+k \cdot m_{22}\\
   (k^2+1) m_{21}+k \cdot m_{11} & (k^2+1)m_{22}+k \cdot m_{12}\end{array}
\right)$. In this case, obviously $|(k^2+1)m_{21}+  k \cdot m_{11}|
> |m_{21}|$ and $|(k^2+1)m_{22}+  k \cdot m_{12}| > |m_{22}|$, so we do not
have a decrease, i.e., this case is moot.

\smallskip

\noindent {\bf Case 1d.} Two subsequent elementary operations are:
 $(row_1 + k \cdot row_2)$, followed by $(row_2 - k \cdot row_1)$. The result of
 applying these two operations to the matrix $M$ is:  $ \left(
\begin{array}{cc} m_{11}+k \cdot m_{21} & m_{12}+k \cdot m_{22}\\
    -(k^2-1)m_{21}- k \cdot m_{11} & -(k^2-1)m_{22}- k \cdot m_{12}\end{array}
\right)$. If $k > \sqrt{2}$, then $|-(k^2-1)m_{21}- k \cdot m_{11}|
> |m_{21}|$ and $|-(k^2-1)m_{22}- k \cdot m_{12}| > |m_{22}|$, so we
do not have a decrease, i.e., this case is moot, too.

\smallskip

\noindent {\bf Case 1e.} Two subsequent elementary operations are:
 $(row_1 - k \cdot row_2)$, followed by $(row_2 + k \cdot row_1)$. The result of
applying these two operations to the matrix $M$ is:  $ \left(
\begin{array}{cc} m_{11}-  k \cdot m_{21} & m_{12}- k \cdot m_{22}\\
    k \cdot m_{11}-(k^2-1) \cdot m_{21}  & k \cdot m_{12}-(k^2-1) \cdot m_{22} \end{array}
\right)$. If the first operation increases $|M|$ (or leaves it
unchanged) and then the second one reduces it, then the second
operation should reduce $|m_{21}|$ or  $|m_{22}|$. Assume, without
loss of generality, that $|m_{21}| \ge |m_{22}|$.

We may assume that $|m_{11}-k \cdot m_{21}| \ge |m_{11}| = -m_{11}$
and $|m_{12}-k \cdot m_{22}| \ge |m_{12}| = m_{12}$ because
otherwise, the complexity of $M$ could be reduced by a single
operation $(row_1 - k \cdot row_2)$.

Now we look at the inequality  $|k \cdot m_{11}-(k^2-1) \cdot
m_{21}| < |m_{21}|=-m_{21}$. Re-write it as follows: $|k
\cdot(m_{11}-k \cdot m_{21}) + m_{21}| < -m_{21}$.  We may assume
that  $m_{11}-k \cdot m_{21}>0$ because otherwise, a single
operation $(row_1 - k \cdot row_2)$ would reduce the complexity of
$M$. We also know that $|m_{11}-k \cdot m_{21}| \ge -m_{11}$ (see
the previous paragraph). Thus, $m_{11}-k \cdot m_{21} \ge -m_{11}$.
This inequality, together with $|k \cdot(m_{11}-k \cdot m_{21}) +
m_{21}| < -m_{21}$, yield   $|-k \cdot m_{11} + m_{21}| < -m_{21}$.
This means  a  single operation $(row_2 - k \cdot row_1)$ would
reduce $|m_{21}|$.
\smallskip

\noindent {\bf Case 1f.}  Two subsequent elementary operations are:
$(row_1 - k \cdot row_2)$, followed by $(row_2 - k \cdot row_1)$.
The result of applying these two operations to the matrix $M$ is:
$\left(
\begin{array}{cc} m_{11}-  k \cdot m_{21} & m_{12}- k \cdot m_{22}\\
   (k^2+1) \cdot m_{21} - k \cdot m_{11}  & (k^2+1) \cdot m_{22} - k \cdot m_{12} \end{array}
\right)$. If the first operation increases $|M|$ (or leaves it
unchanged) and then the second one reduces it, then the second
operation should reduce $|m_{21}|$ or  $|m_{22}|$.

We may assume that $m_{12} - k \cdot m_{22} < 0$ because otherwise,
a single operation $(row_1 - k \cdot row_2)$ would reduce $|m_{12}|$
and therefore also $|m_{11}|$ for the determinant to remain
unchanged in this case. Now look at the element in the bottom right
corner: $(k^2+1) \cdot m_{22} - k \cdot m_{12} = k \cdot (k \cdot
m_{22} -m_{12}) + m_{22}$. Since $k \cdot m_{22} -m_{12} >0$, we
have $|(k^2+1) \cdot m_{22} - k \cdot m_{12}|
> |m_{22}|$, a contradiction.

\smallskip

\noindent {\bf Case 1g.}  Two subsequent elementary operations are:
$(column_1 + k \cdot column_2)$, followed by $(column_2 + k \cdot
column_1)$. The result of applying these two operations to the
matrix $M$ is: $\left(
\begin{array}{cc} m_{11}+  k \cdot m_{12} & (k^2+1) \cdot m_{12}+ k \cdot m_{11}\\
 m_{21} + k \cdot m_{22}  & (k^2+1) \cdot m_{22} + k \cdot m_{21} \end{array}
\right)$. If the first operation increases $|M|$ (or leaves it
unchanged) and then the second one reduces it, then the second
operation should reduce $|m_{12}|$ or  $|m_{22}|$. Let us assume
here that $|m_{22}| \ge |m_{12}|$.

Now look at the element in the bottom right corner: $(k^2+1) \cdot
m_{22}  + k \cdot m_{21} = k \cdot (m_{21} + k \cdot m_{22}) +
m_{22}$. We may assume that $m_{21} + k \cdot m_{22} > 0$ because
otherwise, a single operation $(column_1 + k \cdot column_2)$ would
reduce $|m_{21}|$. In that case, however, we have $|(k^2+1) \cdot
m_{22}  + k \cdot m_{21}| = k \cdot (m_{21} + k \cdot m_{22}) +
m_{22} > m_{22} = |m_{22}|$, a contradiction.
\smallskip

\noindent {\bf Case 1h.}  Two subsequent elementary operations are:
$(column_1 + k \cdot column_2)$, followed by $(column_2 - k \cdot
column_1)$. The result of applying these two operations to the
matrix $M$ is: $\left(
\begin{array}{cc} m_{11}+  k \cdot m_{12} & (1-k^2) \cdot m_{12}- k \cdot m_{11}\\
 m_{21} + k \cdot m_{22}  & (1-k^2) \cdot m_{22} - k \cdot m_{21} \end{array}
\right)$. If the first operation increases $|M|$ (or leaves it
unchanged) and then the second one reduces it, then the second
operation should reduce $|m_{12}|$ or  $|m_{22}|$. Assume here that
$|m_{22}| \ge |m_{12}|$. Then we should have  $|(1-k^2) \cdot m_{22}
- k \cdot m_{21}| < |m_{22}| = m_{22}$.

We may assume that $m_{21} + k \cdot m_{22} > 0$ because otherwise,
a single operation $(column_1 + k \cdot column_2)$ would reduce
$|m_{21}|$ while keeping the element in this position negative. Then
this same operation should reduce $|m_{11}|$, too, for the
determinant to remain unchanged. Also, since the first operation was
supposed to increase $|M|$ (or leave it unchanged), we should have,
in particular, $|m_{21} + k \cdot m_{22}| = m_{21} + k \cdot m_{22}
\ge |m_{21}| = -m_{21}$. This, together with the inequality in the
previous paragraph, gives $|(1-k^2) \cdot m_{22} - k \cdot m_{21}| =
|m_{22} - k \cdot (m_{21} + k \cdot m_{22})| \ge |k \cdot m_{22} +
m_{21}|$. Therefore, we should have $|k \cdot m_{22} + m_{21}|  <
|m_{22}| = m_{22}$. Since $k \ge 2$, this implies $|m_{21}|  >
|m_{22}|$, contradicting the assumption of $m_{22}$ having the
maximum absolute value in the matrix $M$.

\smallskip

\noindent {\bf Case 1i.}  Two subsequent elementary operations are:
$(column_1 - k \cdot column_2)$, followed by $(column_2 + k \cdot
column_1)$. The result of applying these two operations to the
matrix $M$ is: $\left(
\begin{array}{cc} m_{11}-  k \cdot m_{12} & (1-k^2) \cdot m_{12}+ k \cdot m_{11}\\
 m_{21} - k \cdot m_{22}  & (1-k^2) \cdot m_{22} + k \cdot m_{21} \end{array}
\right)$. Since $k \ge 2$, we have $|(1-k^2) \cdot m_{22} + k \cdot
m_{21}| > |m_{22}|$, so this case is moot.
\smallskip

\noindent {\bf Case 1j.}  Two subsequent elementary operations are:
$(column_1 - k \cdot column_2)$, followed by $(column_2 - k \cdot
column_1)$. The result of applying these two operations to the
matrix $M$ is: $\left(
\begin{array}{cc} m_{11}-  k \cdot m_{12} & (1+k^2) \cdot m_{12}- k \cdot m_{11}\\
 m_{21} - k \cdot m_{22}  & (1+k^2) \cdot m_{22} - k \cdot m_{21} \end{array}
\right)$. Since $|(1+k^2) \cdot m_{22} - k \cdot m_{21}| = |m_{22} -
k \cdot (m_{21} - k \cdot m_{22})| > |m_{22}|$, this case is moot,
too.

\medskip

\noindent {\bf Case 2.} Two negative entries are on a diagonal.
Without loss of generality, we assume here that $m_{11}>0, m_{22}>0,
m_{12}<0, m_{21}<0$. Because of the ``symmetry" between row and
column  operations in this case, we can reduce the number of
subcases (compared to Case 1 above) and only consider the following.
\medskip

\noindent {\bf Case 2a.} Two subsequent elementary operations are:
$(row_1 + k \cdot row_2)$, followed by $(column_1 + k \cdot
column_2)$. The result of applying these two operations to the
matrix $M$ is:  $ \left(
\begin{array}{cc} m_{11} + k \cdot m_{21} +k \cdot m_{12} + k^2 m_{22} & m_{12}+k \cdot m_{22}\\
  m_{21}+k \cdot m_{22} & m_{22} \end{array}
\right)$. This case is nontrivial only if the first operation
increases the absolute value of the  element in the top left corner
(or leaves it unchanged), and then the second operation reduces it.
First let us look at the element $m_{21}+k \cdot m_{22}$. If
$m_{21}+k \cdot m_{22} <0$, then a single operation  $(column_1 + k
\cdot column_2)$ would reduce $m_{21}$ while keeping that element
negative. In that case, this operation would reduce $m_{11}$, too,
while keeping it positive because otherwise, the determinant would
change. Thus, a single operation  $(column_1 + k \cdot column_2)$
would reduce the complexity of $M$ in that case. The same argument
shows that $m_{12}+k \cdot m_{22} \ge 0$ because otherwise, a single
operation $(row_1 + k \cdot row_2)$ would reduce the complexity of
$M$.

If $m_{21}+k \cdot m_{22} \ge 0$  and  $m_{12}+k \cdot m_{22} \ge
0$, then $m_{11} + k \cdot m_{21} +k \cdot m_{12} + k^2 m_{22} \ge
0$ for the determinant to be equal to 1. But then the second
operation should also reduce the absolute value of $m_{21}$ for the
determinant to remain unchanged. This means a single operation
$(column_1 + k \cdot column_2)$ would reduce $|m_{21}|$ to begin
with,  and this same operation should also reduce $|m_{11}|$ for the
determinant to remain unchanged, so this single operation would
reduce the complexity of $M$.

The same argument takes care of any of the following pairs of
subsequent elementary operations: $(row_1 \pm k \cdot row_2)$,
followed by $(column_1 \pm k \cdot column_2)$, as well as  $(row_1
\pm k \cdot row_2)$, followed by $(column_2 \pm k \cdot column_1)$.

\smallskip

\noindent {\bf Case 2b.} Two subsequent elementary operations are:
$(row_1 + k \cdot row_2)$, followed by $(row_2 + k \cdot row_1)$.
The result of applying these two operations to the matrix $M$ is:
$\left(
\begin{array}{cc} m_{11}+k \cdot m_{21} & m_{12}+k \cdot m_{22}\\
   (k^2+1) m_{21}+k \cdot m_{11} & (k^2+1)m_{22}+k \cdot m_{12}\end{array}
\right)$. If the first operation increases $|M|$ (or leaves it
unchanged) and then the second one reduces it, then the second
operation should reduce $|m_{21}|$ or  $|m_{22}|$. Assume, without
loss of generality, that $|m_{21}| \ge |m_{22}|$.


Thus, we have $|(k^2+1)m_{21}+  k \cdot m_{11}| < |m_{21}| =
-m_{21}$. At the same time, we may assume that $|m_{11}+k \cdot
m_{21}| \ge |m_{11}| = m_{11}$  because otherwise, a single
operation $(row_1 + k \cdot row_2)$ would reduce the complexity of
$M$.

Then, if $m_{11}+k \cdot m_{21} >0$, then a single operation $(row_1
+ k \cdot row_2)$ would reduce $|m_{11}|$, and therefore also
$|m_{12}|$. Thus, we may assume that $m_{11}+k \cdot m_{21} <0$.

Now let us look at the inequality $|(k^2+1)m_{21}+  k \cdot m_{11}|
< -m_{21}$. We know that $m_{11}+k \cdot m_{21} <0$. Therefore, $k
\cdot m_{11}+k^2 \cdot m_{21} <0$. Now $|(k^2+1)m_{21}+  k \cdot
m_{11} = |k \cdot m_{11}+k^2 \cdot m_{21} + m_{21}| > -m_{21}$. This
contradiction completes Case 2b.

\smallskip

\noindent {\bf Case 2c.} Two subsequent elementary operations are:
 $(row_1 + k \cdot row_2)$, followed by $(row_2 - k \cdot row_1)$. The result of
 applying these two operations to the matrix $M$ is:  $ \left(
\begin{array}{cc} m_{11}+k \cdot m_{21} & m_{12}+k \cdot m_{22}\\
    -(k^2-1)m_{21}- k \cdot m_{11} & -(k^2-1)m_{22}- k \cdot m_{12}\end{array}
\right)$. The analysis here is similar to the previous Case 2b.
First we note that we may assume $m_{11}+k \cdot m_{21} <0$, so  $-k
\cdot m_{11}-k^2 \cdot m_{21} >0$. On the other hand, we may assume
that $|m_{11}+k \cdot m_{21}| = -m_{11}-k \cdot m_{21} \ge m_{11}$
because  otherwise, a single operation $(row_1 + k \cdot row_2)$
would reduce $|m_{11}|$. Therefore, $-k \cdot m_{11} - k^2 \cdot
m_{21} \ge k \cdot m_{11}$.

This, together with the inequality $|-(k^2-1)m_{21}- k \cdot m_{11}|
= |-k \cdot m_{11}-k^2 \cdot m_{21} + m_{21}| < |m_{21}| = -m_{21}$,
implies $|k \cdot m_{11} + m_{21}| < |m_{21}|$, in which case a
single operation $(row_2 + k \cdot row_1)$ would reduce $|m_{21}|$.

\smallskip

\noindent {\bf Case 2d.} Two subsequent elementary operations are:
 $(row_1 - k \cdot row_2)$, followed by $(row_2 + k \cdot row_1)$. The result of
 applying these two operations to the matrix $M$ is:  $ \left(
\begin{array}{cc} m_{11}-k \cdot m_{21} & m_{12}-k \cdot m_{22}\\
    -(k^2-1)m_{21}+ k \cdot m_{11} & -(k^2-1)m_{22}+ k \cdot m_{12}\end{array}
\right)$. If $k^2>2$, here we obviously have $|-(k^2-1)m_{22}+ k
\cdot m_{12}| > |m_{22}|$,  and $|-(k^2-1)m_{21}+ k \cdot m_{11}| >
|m_{21}|$. Thus, this case is moot.
\smallskip

\noindent {\bf Case 2e.} Two subsequent elementary operations are:
 $(row_1 - k \cdot row_2)$, followed by $(row_2 - k \cdot row_1)$.
The result of applying these two operations to the matrix $M$ is:
$\left(
\begin{array}{cc} m_{11}-k \cdot m_{21} & m_{12}-k \cdot m_{22}\\
    (k^2+1)m_{21}- k \cdot m_{11} & (k^2+1)m_{22}- k \cdot m_{12}\end{array}
\right)$. Here again we  obviously have $|(k^2+1)m_{22}- k \cdot
m_{12}| > |m_{22}|$ and $|(k^2+1)m_{21}- k \cdot m_{11}| >
|m_{21}|$, so this case is moot, too.

\medskip

\noindent {\bf Case 3.} There are no negative entries. Because of
the obvious symmetry, it is sufficient to consider the following
cases.
\medskip

\noindent {\bf Case 3a.} Two subsequent elementary operations are
both $(column_1 + k\cdot column_2)$, $(column_2 + k\cdot column_1)$,
$(row_1 + k\cdot row_2)$ or $(row_2 + k\cdot row_1)$. Since all the
entries are positive, the second operation cannot decrease the
complexity in this case.
\smallskip

\noindent {\bf Case 3b.} Two subsequent elementary operations are
$(column_1 + k \cdot column_2)$, followed by $(column_2 - k \cdot
column_1)$, give the matrix\par $\left(
\begin{array}{cc} m_{11}+  k \cdot m_{12} & (1-k^2) \cdot m_{12}- k \cdot m_{11}\\
 m_{21} + k \cdot m_{22}  & (1-k^2) \cdot m_{22} - k \cdot m_{21} \end{array}
\right)$. If the first operation increases $|M|$ (or leaves it
unchanged) and then the second operation reduces it, then the second
operation should reduce, say, $|m_{12}|=m_{12}$ (assuming that
$m_{12}\ge m_{22}$). 

After the first operation we note that $|m_{11}+k\cdot
m_{12}|=m_{11}+k\cdot m_{12}\geq |m_{12}|=m_{12}$. After the second
operation, assuming that complexity was reduced, we have
$|(1-k^2)m_{12}-k\cdot m_{11}|\leq |m_{12}|=m_{12}$, which can be
rewritten as $|m_{12}-k(m_{11}+k\cdot m_{12})|\leq m_{12}$. Since
all the entries are positive, we have $m_{11}+k\cdot m_{12}\geq
m_{12}$, and hence in order for the second operation to reduce
$|m_{12}|$, the following inequality should hold:  $-m_{12}\leq
(1-k^2)m_{12}-k\cdot m_{11}\leq 0$. Subtracting $m_{12}$,
multiplying each part by $-1$ and factoring out $k$ we get
$2m_{12}\geq k(m_{11}+\cdot m_{12})\geq m_{12}$. Divide by $k$:
$\frac{2}{k}\cdot m_{12}\geq m_{11}+k\cdot m_{12}\geq
\frac{1}{k}\cdot m_{12}$. However, $k\geq 2$ implies that
$m_{11}+k\cdot m_{12}\leq m_{12}$, which brings us to a
contradiction.

\smallskip

\noindent {\bf Case 3c.} Two subsequent elementary operations are
$(column_1 - k \cdot column_2)$, followed by $(column_2 + k \cdot
column_1)$. The resulting matrix is \par $\left(
\begin{array}{cc} m_{11}-  k \cdot m_{12} & (1-k^2) \cdot m_{12}+ k \cdot m_{11}\\
 m_{21} - k \cdot m_{22}  & (1-k^2) \cdot m_{22} + k \cdot m_{21} \end{array}
\right)$. If the first operation increases $|M|$ (or leaves it
unchanged) and then the second one reduces it, then the second
operation should reduce, say, $|m_{12}|=m_{12}$ (assuming that
$m_{12}\ge m_{22}$).

Also, after the first operation we can observe that
$m_{11}\leq|m_{11}-k\cdot m_{12}|=$ and $m_{21}\leq |m_{21}-k\cdot
m_{11}|$ because otherwise,  a single operation $(column_1 - k\cdot
column_2)$ would reduce $|m_{11}|$ and $|m_{21}|$. This implies that
$m_{11}-k\cdot m_{12}\leq -m_{11}$ and $m_{21}-k\cdot m_{22}\leq
-m_{21}$.

Consider the inequality $|(1-k^2) m_{12}+k\cdot m_{11}|\leq m_{12}$.
Rewrite it in the following way: $|m_{12}-k\cdot (k\cdot
m_{12}-m_{11})|\leq m_{12}$. Combining this with the previous
inequalities we get $|m_{12}-k\cdot m_{11}|\leq |m_{12}-k\cdot
(k\cdot m_{12}-m_{11})|\leq m_{12}$. Therefore, a single operation
$(column_2 - k\cdot column_1)$ reduces $|m_{12}|$.

\smallskip

\noindent {\bf Case 3d.} Two subsequent elementary operations are
$(column_1 - k \cdot column_2)$ followed by $(column_2 - k \cdot
column_1)$.  The resulting matrix is  \par $\left(
\begin{array}{cc} m_{11}-  k \cdot m_{12} & (1-k^2) \cdot m_{12}- k \cdot m_{11}\\
 m_{21} - k \cdot m_{22}  & (1-k^2) \cdot m_{22} - k \cdot m_{21} \end{array}
\right)$. If the first operation increases $|M|$ (or leaves it
unchanged) and then the second operation reduces it, then the second
one should reduce, say, $|m_{12}|=m_{12}$ (assuming that $m_{12}\ge
m_{22}$).

We may assume that $m_{11}-k\cdot m_{12} \le 0$ because if
$m_{11}-k\cdot m_{12} >0$, then for the determinant to be unchanged
after the first operation, we would have also $m_{21}-k\cdot m_{22}
>0$, but then a single operation
$(column_1 - k\cdot column_2)$ would reduce the complexity of $M$.

Thus, $|m_{11}-k\cdot m_{12}| = k\cdot m_{12} - m_{11}$. If the
first operation $(column_1 - k \cdot column_2)$ did not reduce the
complexity of $M$, then  $m_{11}\leq k\cdot m_{12} - m_{11}$, which
implies that $m_{11}-k\cdot m_{12} \leq -m_{11}$.

Now consider the inequality $|(1-k^2) m_{12}-k\cdot m_{11}| \leq
m_{12}$. After rewriting it we get $|m_{12}-k\cdot (k\cdot
m_{12}-m_{11})|\leq m_{12}$. Combining it with the inequality in the
previous paragraph,  we get $|m_{12}-k\cdot m_{11}|\leq
|m_{12}-k\cdot (k\cdot m_{12}-m_{11})|\leq m_{12}$. Therefore, a
single elementary operation $(column_2 - k\cdot column_1)$ reduces
$|m_{12}|$.

\section{Proof of Theorem \ref{infinite}}

Let $k, m \in \Z, ~k \ge 3, ~m \ge 1$. Denote  $M(k, m)= \left(
\begin{array}{cc} 1-k^2m & k^2m \\ -k^2m & 1+k^2m\end{array}
\right)$. It is straightforward to check that:

\medskip

\noindent {\bf (1)} $M(k, m)$ has determinant 1;

\noindent {\bf (2)} $M(k, m)= M(k, 1)^m$;

\noindent {\bf (3)} No elementary $k$-operation reduces the absolute
value of any of the entries of $M(k, m)$.

\medskip

Since the cyclic group generated by $M(k, 1)$ is infinite, the
result follows from Theorem \ref{peak}.

\section{Proof of Corollary \ref{complexity}}
\label{Corollary}

We assume in this section that $k \ge 3$ because for $k=2$, the
membership problem in the subgroup of $SL_2(\Z)$ generated by $A(k)$
and $B(k)$ is solvable in constant time, see Corollary
\ref{corollary2} in the Introduction.

First of all we check that a given matrix $M$ from the group
$SL_2(\Z)$ has the form $\left(
\begin{array}{cc} 1+k^2n_1 & kn_2 \\ kn_3 & 1+k^2n_4\end{array}
\right)$ for some integers $n_i$.  Then we check that $M$ has at
most one zero entry. If there are more, then $M$ does not belong to
the subgroup in question unless $M$ is the identity matrix. We also
check that $\max |m_{ij}| > 1$. If $\max |m_{ij}| = 1$, then $M$
does not belong to the subgroup in question unless $M$ is the
identity matrix. Indeed, the only nontrivial cases here are $M=A(\pm
1)$ and $M=B(\pm 1)$. Then $M^k=A(k)^{\pm 1}$ or $M^k=B(k)^{\pm 1}$.
This would give a nontrivial relation in the group generated by
$A(k)$ and $B(k)$ contradicting the fact that this group is free.

Now let $\max |m_{ij}| > 1$. If no elementary operation either
reduces $\max |m_{ij}|$ or reduces $\sum_{i,j}|m_{ij}|$ without
increasing $\max |m_{ij}|$, then $M$ does not belong to the subgroup
generated by $A(k)$ and $B(k)$. If there is an elementary operation
that reduces  $\max |m_{ij}|$, then we apply it. For example,
suppose the elementary operation $(row_1 -k \cdot row_2)$ reduces
$|m_{11}|$. The result of this operation is the matrix $\left(
 \begin{array}{cc} m_{11}-km_{21} & m_{12}-km_{22} \\ m_{21} & m_{22} \end{array}
 \right)$. If $|m_{11}|$ here has decreased, then $|m_{12}|$ could not
increase because otherwise, the determinant of the new matrix would
not be equal to 1. Thus, the complexity of the matrix $M$ has been
reduced, and the new matrix belongs to our subgroup if and only if
the matrix $M$ does. Since there are only finitely many numbers of
the form $kn, n \in \Z,$ with bounded absolute value, this process
should terminate either with a non-identity matrix whose complexity
cannot be reduced or with the identity matrix. In the latter case,
the given matrix was in the subgroup generated by $A(k)$ and $B(k)$;
in the former case, it was not.

To estimate the time complexity of this algorithm, we note that each
step of it (i.e., applying a single elementary operation) takes time
$O(\log m)$, where $m$ is the complexity of the matrix this
elementary operation is applied to. This is because if $k$ is an
integer, multiplication by $k$ amounts to $k-1$ additions, and each
addition of integers not exceeding $m$ takes time $O(\log m)$. Since
the complexity of a matrix is reduced at least by 1 at each step of
the algorithm, the total complexity is $O(\sum_{k=1}^n \log k) = O(n
\cdot \log n)$. This completes the proof. $\Box$

\smallskip

As for generic-case complexity of this algorithm (cf. Problem
\ref{generic} in our Section \ref{results}), we note that, speaking
very informally, a ``random" product of $A(k)^{\pm 1}$ and
$B(k)^{\pm 1}$ is ``close" to a product where $A(k)^{\pm 1}$ and
$B(k)^{\pm 1}$ alternate, in which case the complexity of the
product matrix grows exponentially in the number of factors (see
e.g. \cite[Proposition 1]{BSV}), so the number of summands in the
sum that appears in the proof of Corollary \ref{complexity} will be
logarithmic in $n$, and therefore generic-case complexity of the
algorithm should be $O(\log ^2 n)$ in case the answer is ``yes"
(i.e., an input matrix belongs to the subgroup generated by $A(k)$
and $B(k)$). Of course, a ``random" matrix from $SL_2(\Z)$ will {\it
not} belong to the subgroup generated by $A(k)$ and $B(k)$ with
overwhelming probability. This is because if $k \ge 3$, this
subgroup has infinite index in $SL_2(\Z)$. It is, however, not clear
how fast (generically) our algorithm will detect that; specifically,
whether it will happen in sublinear time or not.

Note that, unlike the algorithms with low generic-case complexity
considered in \cite{KMSS}, this algorithm has a good chance to have
low generic-case complexity giving both ``yes" and ``no" answers.

Finally, we note that generic-case complexity depends on how one
defines the {\it asymptotic density} of a subset of inputs in the
set of all possible inputs. This, in turn, depends on how one
defines the complexity of an input. In \cite{KMSS}, complexity of an
element of a group $G$ was defined as the minimum word length of
this element with respect to a fixed generating set of $G$. In our
situation, where inputs are matrices over $\Z$, it is probably more
natural to define complexity $|M|$ of a matrix $M$ as the sum of the absolute values of
the entries of $M$, like we did in this paper. Yet another natural way is to use {\it Kolmogorov complexity}, i.e., speaking informally, the minimum possible size of a description of $M$. Since Kolmogorov complexity of an integer $n$ is equivalent to $\log n$, we see that for a matrix $M \in SL_2(\Z)$, Kolmogorov complexity is equivalent
to $\log |M|$, for $|M| = \sum |m_{ij}|$, as defined in this paper. This is not the case though if $M \in SL_2(\Q)$ since for Kolmogorov complexity of a rational number, complexity of both the numerator and denominator matters.

\bigskip

\noindent {\it Acknowledgement.} We are grateful to Norbert A'Campo,
Ilya Kapovich, and Linda Keen for helpful comments.

\end{document}